\documentclass[12pt]{amsproc}

\usepackage[margin = 1in]{geometry}
\usepackage{amsmath,amsfonts,amssymb,amsthm}
\usepackage[usenames,dvipsnames]{xcolor}
\usepackage[osf]{mathpazo}
\usepackage{graphicx}

\usepackage{hyperref}
\hypersetup{
    colorlinks,	%
    citecolor=blue,%
    filecolor=black,%
    linkcolor=red,%
    urlcolor=gray
}

	\newtheorem*{theorem*}{Theorem}

\theoremstyle{definition}

\theoremstyle{remark}
	
	\newtheorem*{remark*}{Remark}

 { \begin{list}%
         {$\bullet$}%
         {\setlength{\labelwidth}{20pt}%
          \setlength{\leftmargin}{25pt}%
          \setlength{\topsep}{0pt}
          \setlength{\itemsep}{1.5ex}
          \setlength{\parsep}{0pt}}}%
 { \end{list} }

\newcommand{\excise}[1]{}

\newcommand{\R}{\mathbb{R}}

\newcommand{\N}{\mathbb{N}}

\begin{document}

\title{Discussion of `Event History and Topological Data Analysis'}

\author{Peter Bubenik}
\address{Department of Mathematics, University of Florida}
\email{peter.bubenik@ufl.edu}



\maketitle
 


\begin{abstract}
Garside et al.~\cite{Garside:2021} use event history methods to analyze topological data. We provide additional background on persistent homology to contrast the hazard estimators used by Garside et al.\ with traditional approaches in topological data analysis. In particular, the former is a local method, which has advantages and disadvantages, while homology is a global. We also provide more background on persistence landscapes and show how a more complete use of this statistic improves its performance.
\end{abstract}

\section{Critical points and persistent homology} \label{sec:ph}

Consider a function $f: \R^d \to \R$. For simplicity, let us restrict to the case $d=2$. We are interested in the topology or `shape' of the \emph{lower excursion sets} or \emph{sublevel sets} $\{x \in \R^2 \ | \ f(x) \leq a \}$ for $a \in \R$. Consider, for example, the negative of a kernel density estimator. In particular, we compute the \emph{homology} of these sublevel sets with binary coefficients. The result is two vector spaces $H_0$ and $H_1$, called homology in {degree} $0$ and homology in degree $1$, respectively. The dimension of these vector spaces equals the \emph{Betti numbers}, $\beta_0$ and $\beta_1$. These topological invariants enumerate the number of connected components and the number of `holes', respectively, of the sublevel sets.

In computational settings, one works with a discretization of the domain $\R^2$ given by either triangles or rectangles together with function values on the vertices, edges, and faces. Often the function values are only given on the vertices, which are arranged in a rectangular grid and the values on the edges and faces are set equal to the maximum of the values on their boundary vertices.

If we increase the value of the level $a$ then the homology of the sublevel sets changes at the critical values of the function $f$. At the value of local minima, $\beta_0$ increases by one. At value of local maxima, $\beta_1$ decreases by one. At the value of saddle points, either $\beta_0$ decreases by one or $\beta_1$ increases by one. The local minima, saddle points, and local maxima are said to be \emph{critical points} of \emph{index} $0$, $1$, and $2$, respectively. 

For $a \leq b$ the sublevel set given by $a$ is contained in the sublevel set given by $b$. This inclusion of subspaces of $\R^2$ induces a linear map between the corresponding homology vector spaces. The collection of vector spaces and linear maps (for all $a \leq b$) is called a \emph{persistence module}.
\emph{Persistent homology} provides a complete description of the persistence module by determining a canonical pairing between critical values. This pairing pairs the critical value for which a homological feature appears, called \emph{birth} with the critical value at which that homological feature disappears, called \emph{death}. In particular, we have pairs in degree $0$ for which a connected component is born at a local minimum and is merged with another component at a saddle point and we have pairs in degree $1$ for which a hole is born at a saddle point and which dies at a local maximum. The only remaining critical points are the global minimum and the global maximum. These are either paired (extended persistence) or omitted (reduced homology).
This pairing of critical values by persistent homology is visualized by either a collection of intervals, called a \emph{barcode}, or by plotting the paired critical values as points in the plane, called a \emph{persistence diagram}.

\section{Local versus global} \label{sec:local-global}

The critical values are locally determined. In the smooth, nondegenerate case, the critical points are the points at which the gradient vanishes, and their index is given by the number of negative values of the Hessian. In the discrete case, assuming vertices have distinct values, critical points and their index are determined by comparing neighboring values.
In contrast, one cannot determine whether or not a saddle point increases $\beta_1$
or decreases $\beta_0$
using local information.
In particular, one cannot obtain persistent homology from knowledge of the critical values and their indices.
In the other direction, given persistent homology one may obtain the critical values and their indices.
That is, persistent homology is a strictly more informative statistic than critical values and their indices.
However, critical values and their indices may be quickly computed, while computing persistent homology is comparatively slow.

Therefore, a statistician should choose between the less expensive local statistic -- critical values and their indices -- or the more expensive global statistic -- persistent homology. The appropriate choice depends on the application. For situations in which there is no global structure of interest, a local statistic may be a good choice, while situations in which global `shape' is important call for a global statistic. For an example of the former, consider a stationary isotropic random field. For an example of the latter, consider the reconstruction of neuronal networks in medical images~\cite{Hu:2021}.
The Nelson-Aalen estimator in the manuscript under discussion is a local topological statistic.

\section{Persistence landscapes and their use} \label{sec:pl}

\begin{figure} 
  \centering
  \includegraphics[width=0.3\linewidth]{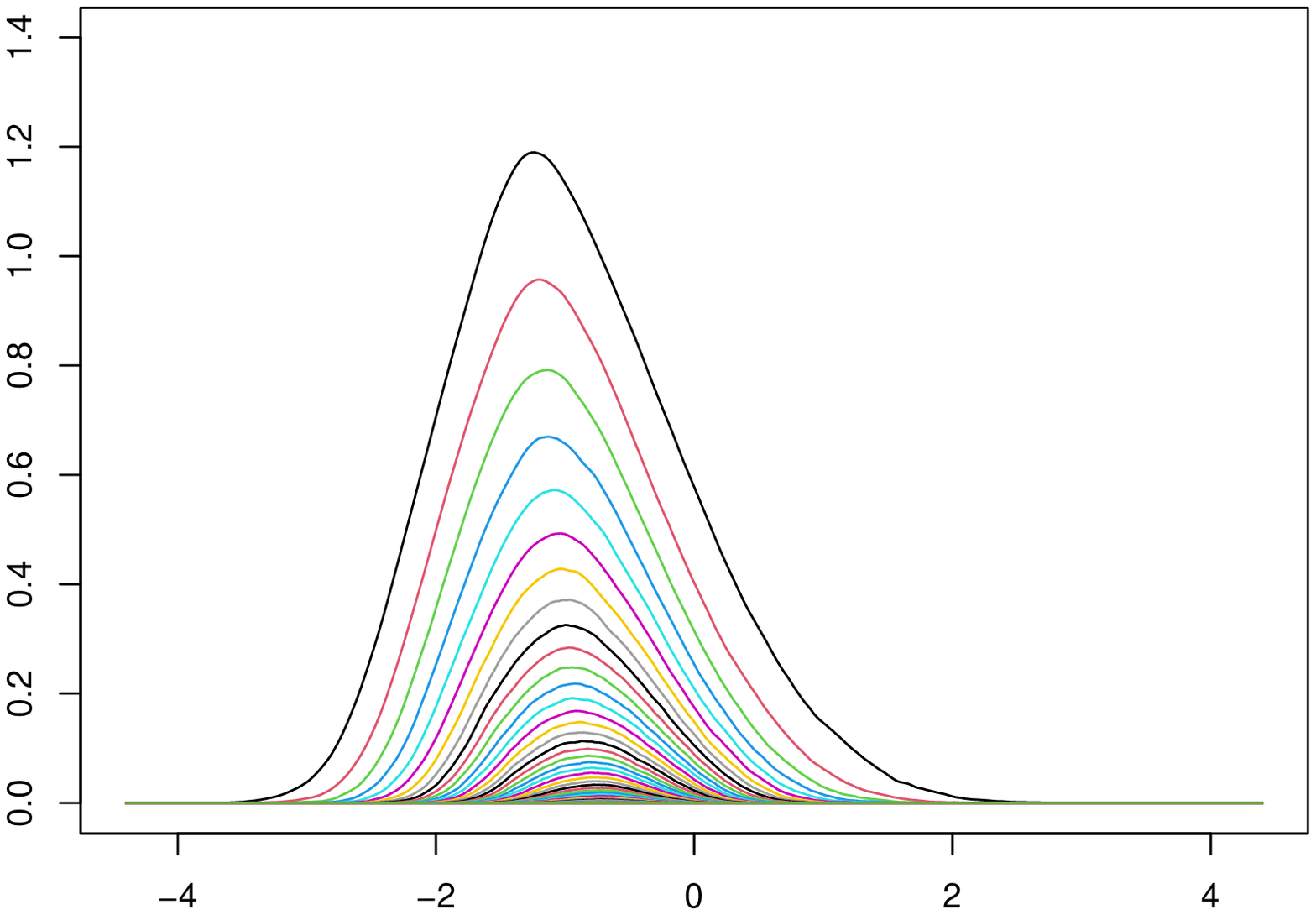} \quad
  \includegraphics[width=0.3\linewidth]{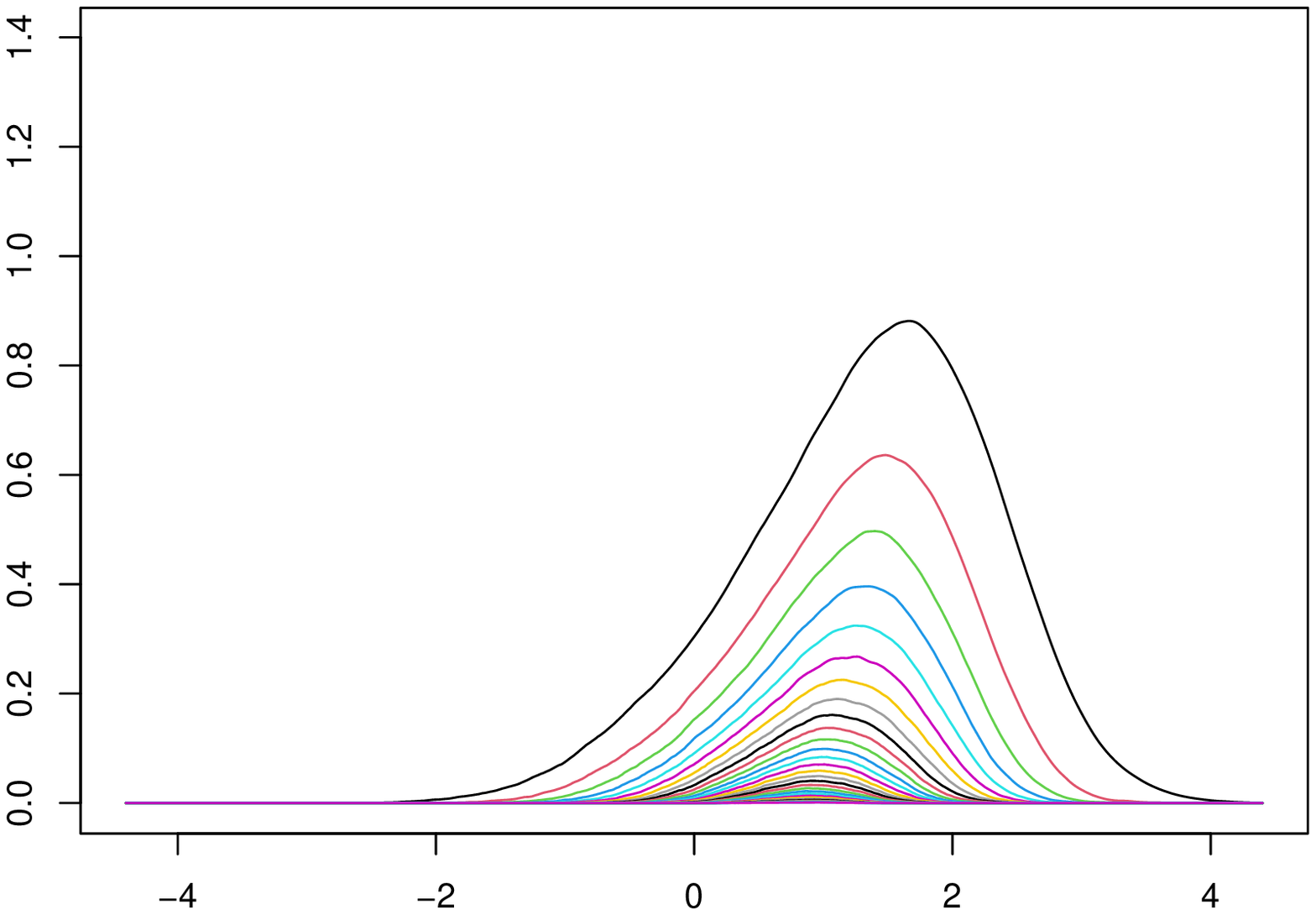}\\
  \includegraphics[width=0.3\linewidth]{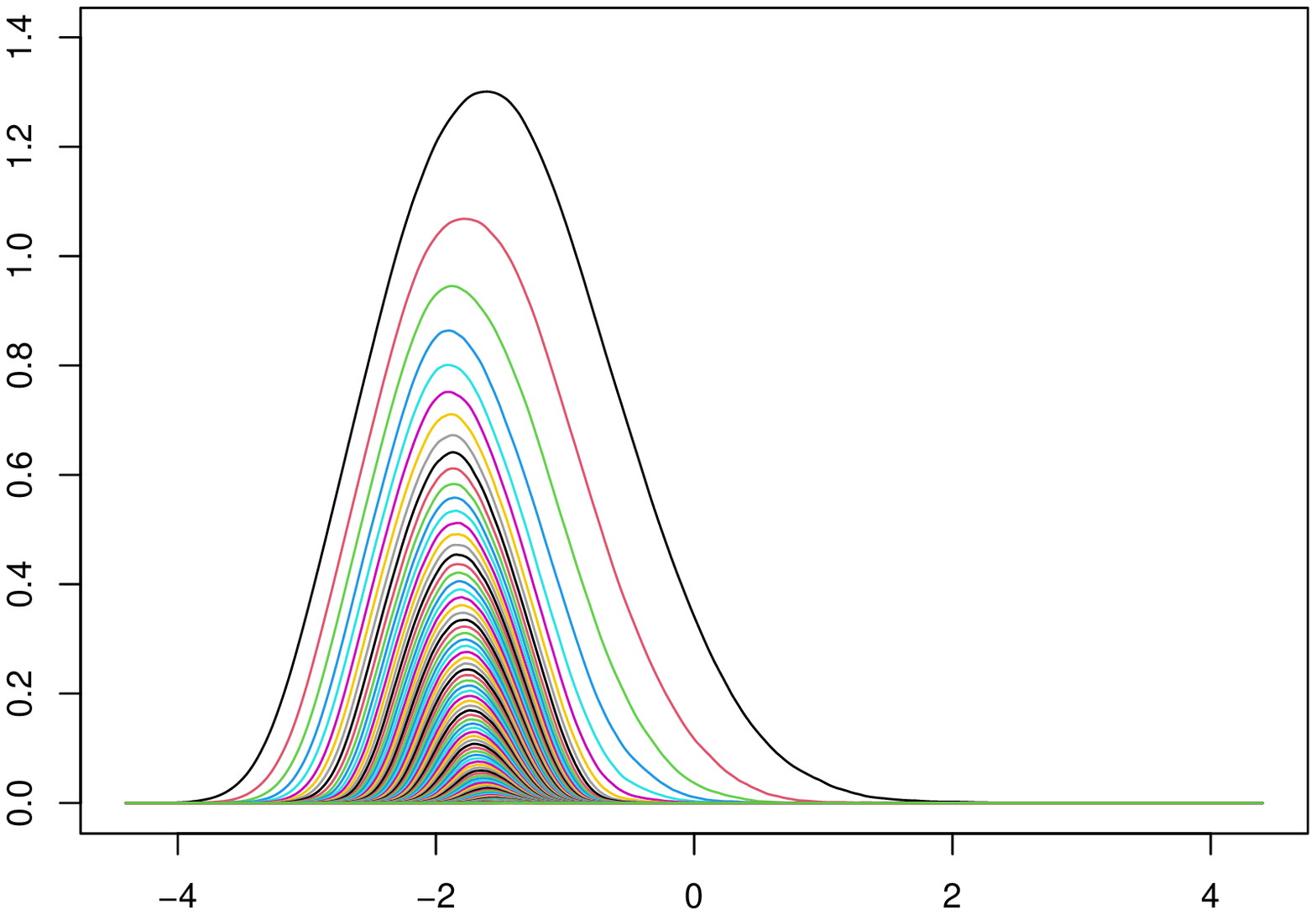} \quad
  \includegraphics[width=0.3\linewidth]{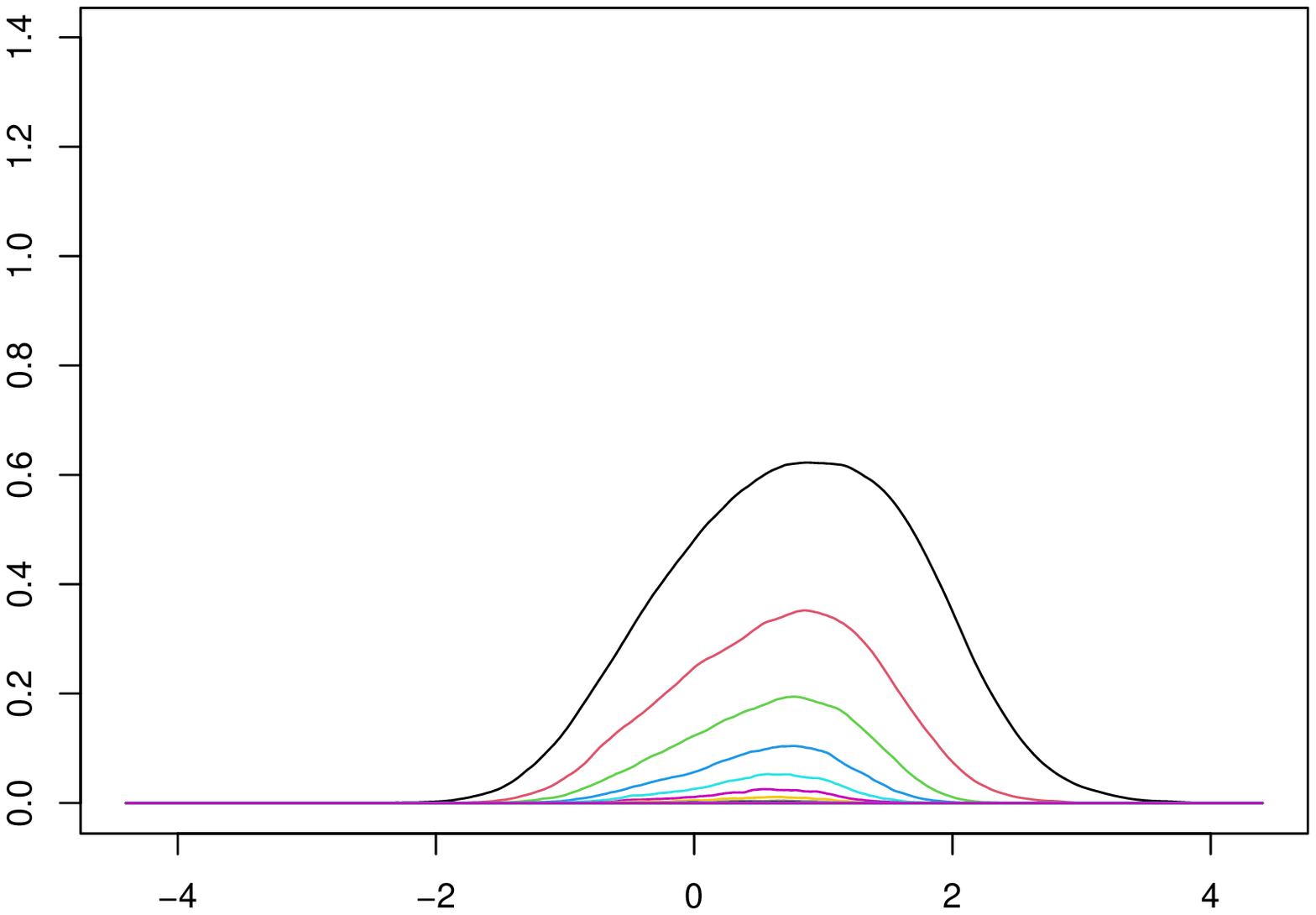}\\
  \includegraphics[width=0.3\linewidth]{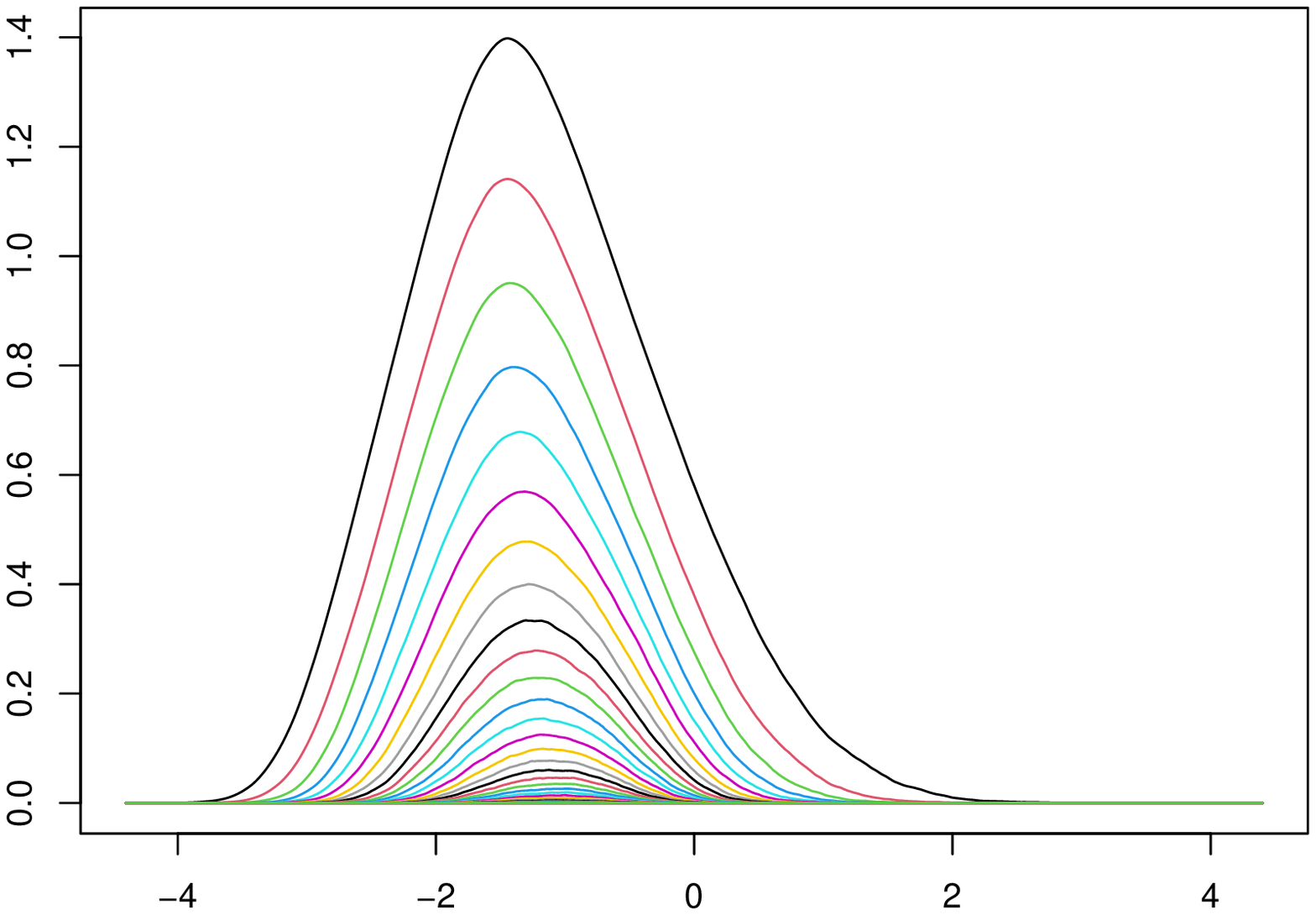} \quad
  \includegraphics[width=0.3\linewidth]{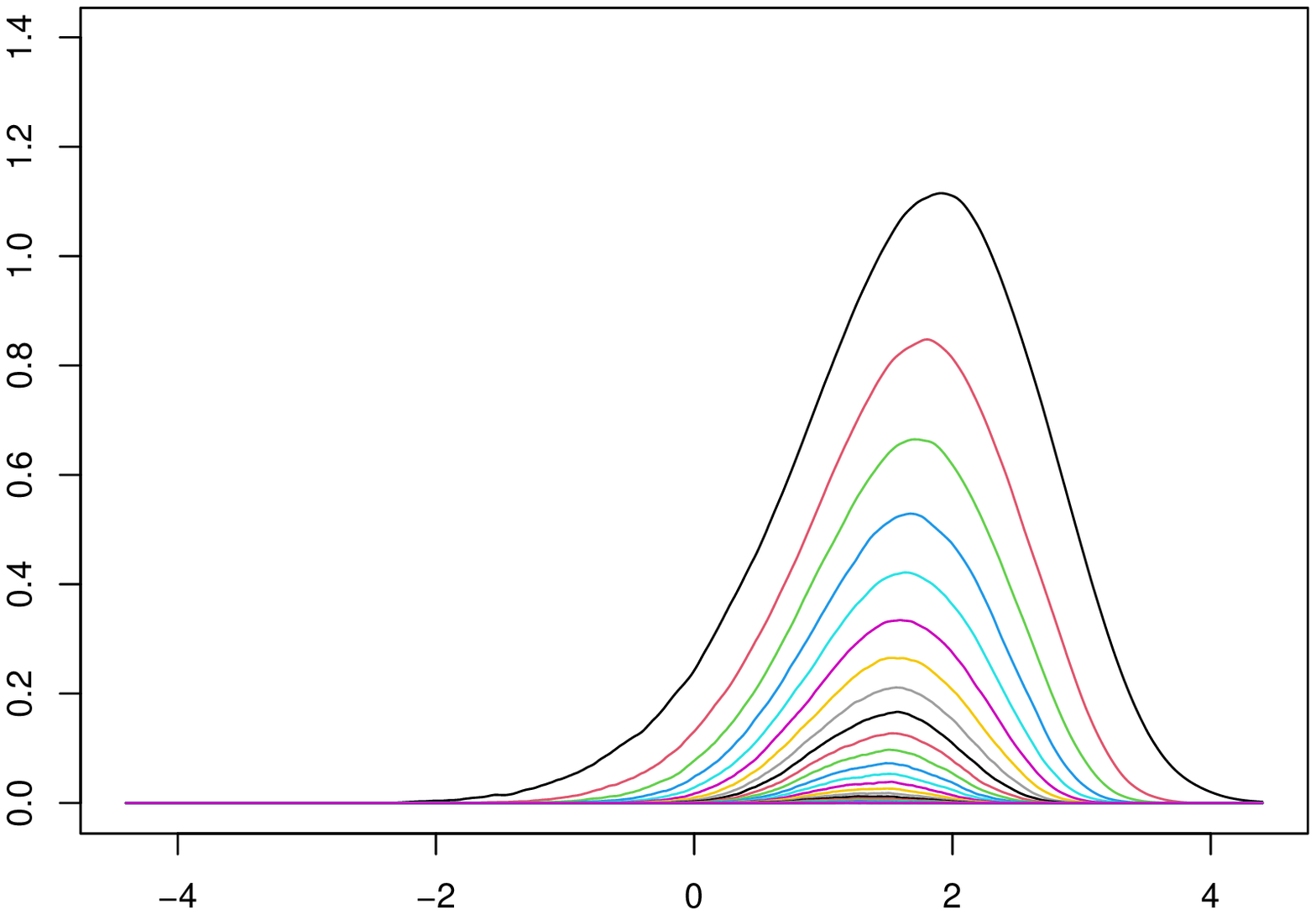}
  \caption{The average persistence landscape in degree $0$ using reduced homology (left) and degree $1$ (right) of 1000 samples of M1 (top row), M2 (middle row), and M3 (bottom row). The $x$ axis gives the function values and the $y$ axis gives the persistence of the homological features of given multiplicity at those values.}
  \label{fig:1}
\end{figure}

\begin{figure}
  \centering
  \includegraphics[width=0.3\linewidth]{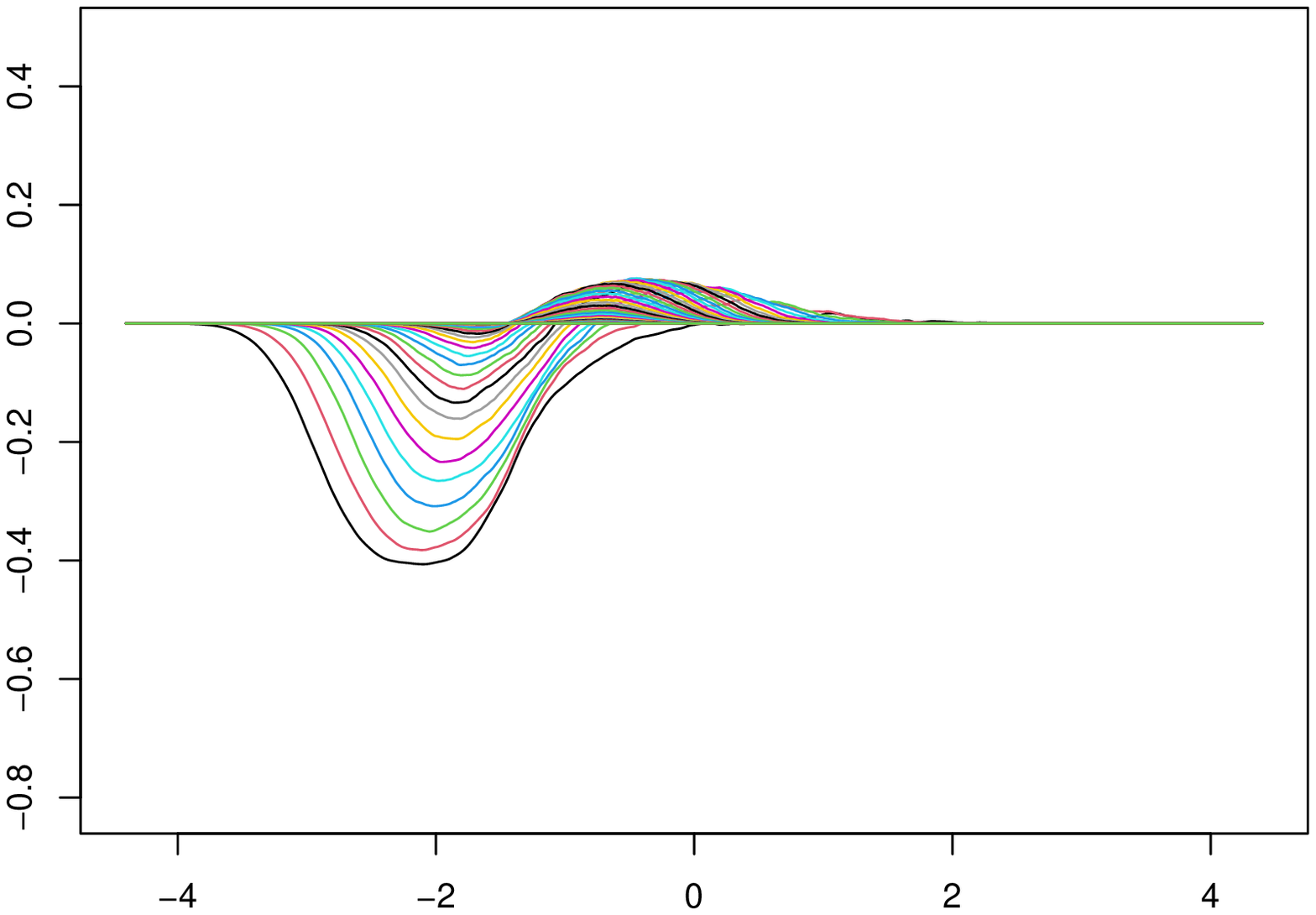} \quad
  \includegraphics[width=0.3\linewidth]{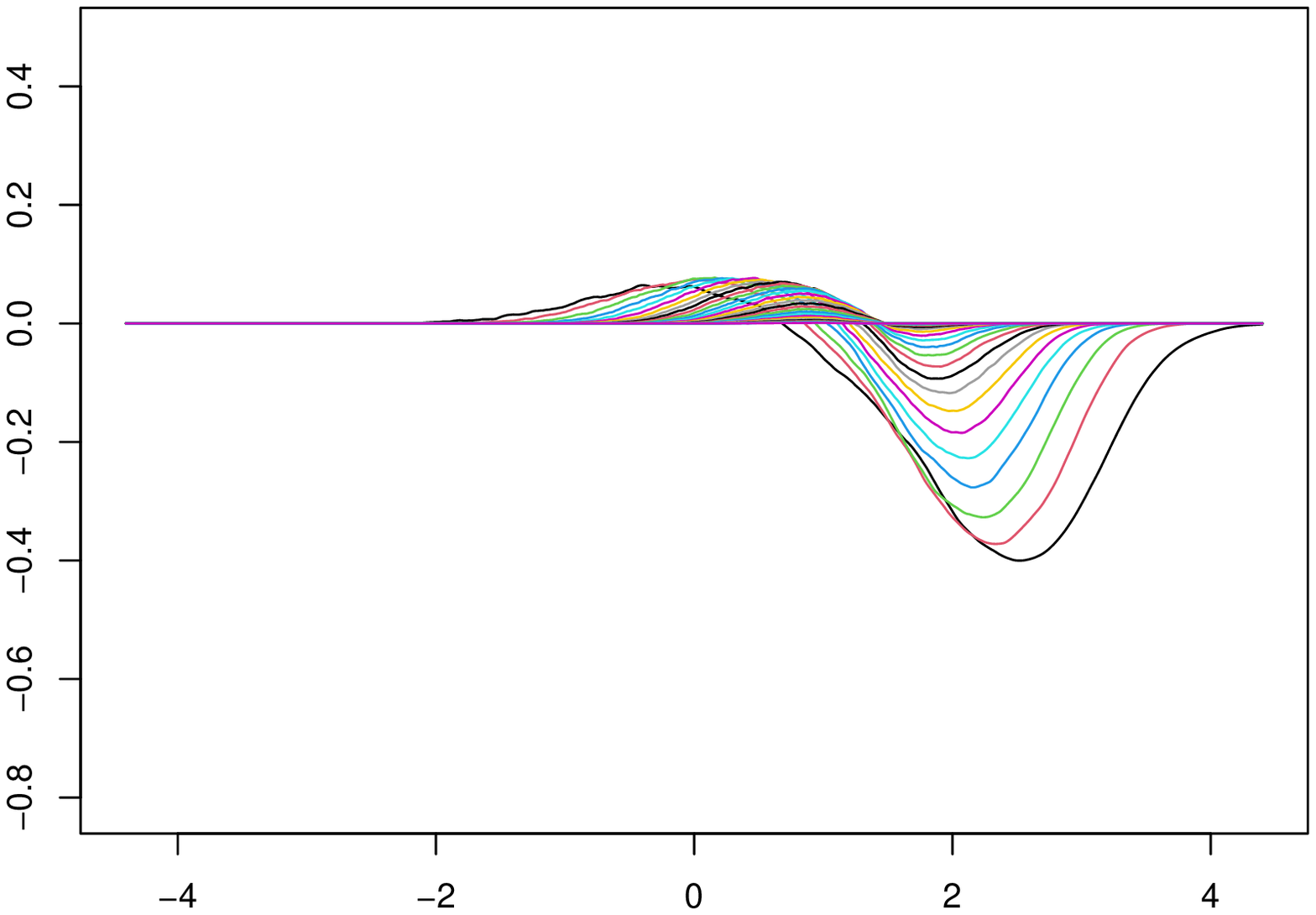} \quad
  \caption{The difference in average persistence landscapes in degree $0$ (left) and degree $1$ (right) between M1 and M3.}
  \label{fig:2}
\end{figure}

Persistence landscapes may be used to provide a vector encoding of persistent homology.
Persistence landscapes are a sequence $(\lambda_k)_{k \geq 1}$ of functions $\lambda_k:\R \to \R$ or, equivalently, a single function $\lambda: \N \times \R \to R$, where $\lambda(k,t) = \lambda_k(t)$.
It is a consequence of their definition that for $k \leq \ell$ and for all $t$, $\lambda_k(t) \geq \lambda_{\ell}(t) \geq 0$ and given a finite amount of data there is a maximum $M$ such that $\lambda_M(t) > 0$ for some $t$.
See Figure~\ref{fig:1} for the average persistence landscapes for the Gaussian random fields considered by Garside et al.~\cite{Garside:2021} The difference between the two most similar cases is given in Figure~\ref{fig:2}.

To obtain a high-dimensional vector representation of persistent homology, we discretize the persistence landscape. In particular, choose values $t_0 < t_1 < \cdots < t_N$ and a choose a maximum \emph{depth} $K$.
Let $\lambda^{(0)}$ and $\lambda^{(1)}$ denote the persistence landscapes for homology in degree $0$ and $1$, respectively.
Define $\Lambda \in \R^{2(N+1)K}$ by
\begin{equation} \label{eq:1}
  \Lambda = (\lambda^{(0)}_1(t_0),\ldots,\lambda^{(0)}_1(t_N),\lambda^{(0)}_2(t_0),\ldots,\lambda^{(0)}_K(t_N),\lambda^{(1)}_1(t_0),\ldots,\lambda^{(1)}_K(t_N)).
\end{equation}
When using persistence landscapes for statistical analysis, the main computational cost is the cost of computing persistent homology. Restricting to a vector $(\lambda^{(0)}_k(t_0),\ldots,\lambda^{(0)}_k(t_0))$ for some fixed $k$ as is done by Garside et al.~\cite{Garside:2021} loses information and does not provide a significant overall savings in computational cost.
Repeating the computations in Garside et al.\ using the vector in \eqref{eq:1}, we have the following updates to their Table 2 in our Table~\ref{tab:1}.
With this update, the Nelson-Aalen estimator and persistence landscape provide  comparable excellent results for M1 v M2 and M2 v M3. The results for the persistence landscape for M1 v M3 are improved but still not as good as those of the Nelson-Aalen estimator. It would be good to understand reason for this difference.

We remark that when using the R package \verb+e1071+ with persistence landscapes, one should use the parameters \verb+scale=FALSE+ and \verb+kernel="linear"+ to use the inner product appropriate for $L^2(\N \times \R)$. Also, since \verb+e1071+ uses sparse matrices internally and persistence landscapes are sparse it is recommended to use a sparse matrix encoding of the persistence landscape.

\begin{table}
  \centering
  \begin{tabular}{ccccccc}
    Mat\'ern & \multicolumn{2}{c}{M1 v M2} & \multicolumn{2}{c}{M1 v M3} & \multicolumn{2}{c}{M2 v M3} \\
    $(\eta,\nu)$ & Acc. & Cal. & Acc. & Cal. & Acc. & Cal.\\
    (5,1) &  100.0 & 97.4 & 93.4 & 88.1 & 100.0 & 97.2\\
    (10,1) &  100.0 & 96.3 & 83.1 & 73.3 & 98.8 & 94.9\\
    (5,2) &  100.0 & 97.4 & 87.9 & 80.9 & 100.0 & 97.1
  \end{tabular}
  \caption{Classification of random fields using persistence landscapes. Accuracy (Acc.) measures the proportion of random fields that were correctly allocated and calibration (Cal.) gives the mean estimate class probabilities for the correct values.}
   \label{tab:1}
\end{table}

 
\bibliographystyle{plain}
\bibliography{biometrika}    

\end{document}